\documentclass[12pt]{amsart}
\usepackage{latexsym, amssymb, amsmath}
\usepackage{eucal}

    \newtheorem{rema}{Remark}[section]
    \newtheorem{propo}[rema]{Proposition}

   \newtheorem{theo}[rema]{Theorem}

   \newtheorem{defi}[rema]{Definition}
    \newtheorem{lemma}[rema]{Lemma}
    \newtheorem{corol}[rema]{Corollary}
     \newtheorem{exam}[rema]{Example}
     \newtheorem{rmk}[rema]{Remark}
        \newcommand{\nno}{\nonumber}

        \newcommand{\p}{\partial}

 \newcommand{\pf}{{\it Proof:}\hspace{2ex}}
 
 \newcommand{\epfv}{\hspace{1em}$\Box$\vspace{1em}}

\newcommand{\bC}{{\mathbb C}}
\newcommand{\bZ}{{\mathbb Z}}
\newcommand{\bQ}{{\mathbb Q}}

\newcommand{\bN}{{\mathbb N}}
\newcommand{\bT}{{\mathbb T}}

\newcommand{\cD}{\mathcal D}
\newcommand{\cP}{\mathcal P}
\newcommand{\cC}{\mathcal C}
\newcommand{\cS}{\mathcal S}

\newcommand{\BQ}{\begin{eqnarray}}
\newcommand{\EQ}{\end{eqnarray}}
\newcommand{\BQn}{\begin{eqnarray*}}
\newcommand{\EQn}{\end{eqnarray*}}

\begin{document}

\title[D-log and formal flow for analytic isomorphisms]
{D-log and formal flow for analytic isomorphisms of n-space}
    \author{David Wright and Wenhua Zhao}

    \bibliographystyle{alpha}
\begin{abstract}  Given a formal map $F=(F_1\ldots,F_n)$ of the form $z+\text{higher}$ 
order terms, we give tree expansion formulas and associated algorithms for the D-Log of 
$F$ and the formal flow $F_t$.  The coefficients which appear in these formulas can be 
viewed as certain generalizations of the Bernoulli numbers and the Bernoulli polynomials.
Moreover the coefficient polynomials in the formal flow formula coincide with the 
strict order polynomials in combinatorics for the partially ordered sets induced 
by trees.  Applications of these formulas to the Jacobian Conjecture are discussed.
\end{abstract}


\renewcommand{\subjclassname}{\textup{2000} Mathematics Subject Classification}
\copyrightinfo{2000}{American Mathematical Sociaty}



\maketitle

\tableofcontents

\renewcommand{\theequation}{\thesection.\arabic{equation}}
\renewcommand{\therema}{\thesection.\arabic{rema}}
\setcounter{equation}{0}
\setcounter{rema}{0}

\section{Introduction}

This work began as an effort to link and extend the results of
\cite{W} and \cite{Z}, placing them in a common framework.  Both of
these papers deal with the formal inverse $F^{-1}$ of a system of
power series $F=(F_1,\ldots,F_n)$; both give formulas for $F^{-1}$ in
terms of $F$, the former being a tree formula, the latter an
exponential formula.  This quest has led to a host of interesting
connections, algorithms, formulas, and relationships with
combinatorics, Bernoulli numbers, and Bernoulli polynomials.

The former paper deals with tree formulas as they apply to formal
inverse, a thread which is also the main thrust of \cite{BCW},
\cite{W1}, \cite{W}, and \cite{C}.  It has combinatoric connections
with generating functions and enumeration techniques for trees.  The
general goal of power series inversion (sometimed called "reversion", 
perhaps to distinguish functional inverse from multiplicative inverse) 
is as follows.  Let $F=(F_{1},\ldots,F_{n})$
with $F_{i}\in\bC[[z_{1},\ldots,z_{n}]]$ for each $i$ and
$F_{i}=z_{1}+\text{terms of degree $\ge2$}$.  One seeks formulas for
the unique $G_{1},\ldots,G_{n}\in\bC[[z_{1},\ldots,z_{n}]]$ for which
$G_{i}(F)=z_{i}$, for $i=1,\ldots,n$.  Perhaps the first of these was
the Lagrange Inversion Formula (see \cite{St2}, Chapter 5), which
dealt with the case $n=1$, and which was generalized (under a certain
restrictive hypothesis) to all $n$ in the by I. J. Good \cite{G} in
1960.  Good then uses his formula for problems of enumerating certain
trees.  In fact, Good's formula had been discovered and published by
Jacobi in 1830 \cite{J}.  Another paper which appeared in 1960 was
that of G. N. Raney \cite{R}, who also related 
formal inverse to trees.  Raney's work was generalized in \cite{C},
which also utilized the work of Jacobi.  A general inversion formula
was given by Abhyankar and Gurjar in 1974 \cite{A}, and this is the
source from which the tree formula of \cite{BCW} was derived, with the
hope of applying it to the Jacobian Conjecture.  Other treatments of
the subject of inversion are \cite {HS}, \cite{Ge}, and \cite{Jo}.

The tree formula of \cite{BCW} expresses the formal inverse $F^{-1}$
as an infinite $\bQ$-linear combination of certain power series
$\cP_T\in \bC[[z]]$, which are constructed using finite rooted trees
$T$.  This construction will be reviewed in \S2 and a new (and quick) proof of the
inversion formula, using the tools developed in this paper,
will be presented in \S5 (Theorem \ref{Theo-5.1}).

Amongst the results of the latter paper is the realization of the $F$
by an expression $F=\text{exp}(A)\cdot z$, where $A=A(z)$, called the
{\it D-Log} of $F$, is a differential operator uniquely determined by
$F$ and yielding the formal inverse as $F^{-1}=\text{exp}(-A)\cdot z$.
Furthermore, the {\it formal flow} $F_t=\text{exp}(tA)\cdot z$
encodes all powers $F^{[n]}$ with $n\in \bZ$ of the formal map $F$.  The
D-Log and the formal flow will be reviewed in \S3.

A primary goal was to show that the D-Log $A$ can also be expressed a
$\bQ$-linear combination of the power series $\cP_T$.  This goal was
attained, yielding a tree formula for the D-Log.  Moreover, we
discovered that the rational coefficients $\phi_T$ of this expression
can be generated by an elegant recurrence relationship and possess
some intriguing combinatorial properties.  For example, the Bernoulli
numbers appear amongst these coefficients.

This situation is placed in a larger context which incorporates formal inverse
by considering the formal
flow $F_t=\text{exp}(tA)\cdot z$, where $t$ is an indeterminate.  For
$n\in\bN$, setting $t=n$ gives the $n$-fold composition
$F\circ\cdots\circ F$, and setting $t=-n$ gives the $n$-fold
composition $F^{-1}\circ\cdots\circ F^{-1}$.  The system $F_t$ can be
written as a $\bQ[t]$-linear combination of the power series $\cP_T$,
producing for each rooted tree $T$ a polynomial $\psi_T(t)$ having
$\phi_T$ as the coefficient of $t$.  Among these polynomials are the
binomial polynomials $\binom {t}{m}$, for all positive integers $m$.
We give an algorithm for calculating $\psi_T(t)$ using the difference
operator $\Delta$.  This formula is used to establish the relationship
of certain $\psi_T(t)$s with the Bernoulli polynomials $B_m(t)$ via an
integration formula.  It is shown that $\psi_T(t)$ also provides an
interesting combinatorial connection: It coincides with the strict
order polynomial $\bar\Omega(P,t)$ (see \cite{St1}, Chapter 4) for
$P=T$, which, for $t=m\in\bZ^+$, counts the number of strict order
preserving maps from any partially ordered set $P$ to the totally
ordered set with $m$ elements.

We would like to thank Professor John Shareshian for informing us of 
Theorem \ref{Omega}, and also Professor Steve Krantz for a helpful 
conversation on flow of analytic maps.

\section {Tree Operations}\label{defs}

\subsection{Notation} By a {\it rooted tree} we mean a finite
1-connected graph with one vertex designated as its {\it root.} The
1-connectivity provides the notion of {\it distance} between two
vertices, which is defined as the number of edges in the unique
geodesic connecting the two.  The {\it height} of a tree is defined to
be the maximum distance of any vertex from the root.  In a rooted tree
there are natural ancestral relations between vertices.  We say a
vertex $w$ is a child of vertex $v$ if the two are connected by an
edge and $w$ lies further from the root than $v$.  In the same
situation, we say $w$ is the {\it parent} of $v$.  Note that a vertex
my have several children, but only one parent.  The root is the only
vertex with no parent.  A vertex is called a {\it leaf}\/ if it has no
children.  When we speak of isomorphisms between rooted trees, we will
always mean root-preserving isomorphisms.

With these notions in mind, we establish the following notation.
\begin{enumerate}
\item We let $\mathbb T$ be the set of isomorphism classes of all
rooted trees and, for $m\ge1$ an integer, we let $\bT_m$
the set of isomorphism classes of all rooted trees with $m$ vertices.
The latter is a finite set.
\item For any rooted tree $T$, we set the following notation:
\begin{itemize}
\item $\text{rt}_T$ denotes the root vertex of $T$.
\item $E(T)$ denotes the set of edges of $T$.
\item $V(T)$ denotes the set of vertices of $T$.
\item $L(T)$ denotes the set of leaves of $T$.
\item $v(T)$ (resp. $l(T)$) denotes the
 number of the elements of $V(T)$ (resp. $L(T)$).
\item $h(T)$ denotes the height of $T$.
\item $\alpha_T$ denotes the number of the elements of the automorphism
group $\mbox{Aut}(T)$.
\item For $v\in V(T)$ we denote by $\alpha_{T,v}$ the size
of the stabilizer of $v$ in $\mbox{Aut}(T)$.  Similarly, for $e\in
E(T)$, we write $\alpha_{T,e}$ the size of the stabilizer of $e$ in
$\mbox{Aut}(T)$.
\item  For $e\in E(T)$ we denote by $v_e$ and $v'_e$ the two (distinct)
vertices which are connected by $e$, with $v_e$ being the one closest to

the root.
\item For $v\in V(T)$ we denote by $v^+$ the set of vertices which are
children of $v$.
\item For $v\in V(T)$ we define the {\it height} of $v$ to be the
number of edges in the (unique) geodesic connecting $v$ to
$\text{rt}_{T}$.
The {\it height} of $T$ is defined to be the maximum of the heights of
its vertices.
\item For $v_{1},\ldots,v_{r}\in V(T)$, we write $T\backslash
\{v_{1},\ldots,v_{r}\}$ for the graph obtained by deleting each of
these vertices and all edges adjacent to these vertices.
\end{itemize}
\item A {\it rooted subtree} of a rooted tree $T$ is defined as a
connected subgraph of $T$ containing $\text{rt}_{T}$, with
$\text{rt}_{T'}=\text{rt}_{T}$.  In this case we write $T'\le T$.  If
$T'\ne T$ we write $T'<T$.  If $T'<T$, we write $T\backslash T'$ for
the graph obtained by deleting all vertices of $T'$ and all edges
adjacent to its vertices.
\item For any $k\geq 1$, we denote by
$C_k$ the rooted tree of height $k-1$ having $k$ vertices, and by
$S_k$ the rooted tree of height $1$ having $k$ leaves.  We also set
$S_0=\circ$, the rooted tree with one vertex.  We refer to the trees
$C_{k}$ as {\it chains} and the $S_{k}$ as {\it shrubs.}
\end{enumerate}

\subsection{Power Series Given by a Rooted Tree} Let
$\bC[[z_1,\ldots,z_n]]=\bC[[z]]$ denote the ring of formal power
series in $n$ indeterminates $z_1,\ldots,z_n$ over the complex numbers\footnote{In this paper $\bC$ can always be replaced by any $\bQ$-algebra.} $\bC$.  For $i=1,\ldots,n$ we will write $D_i$ for the differential
operator $\frac{\partial}{\partial z_i}$.  The operators
$D_1,\ldots,D_n$ are commuting derivations acting on the ring
$\bC[[z]]$.

Given a vector of power series $F=(F_1,\ldots,F_n)\in\bC[[z]]^n$, we
write $F_i=z_i+H_i$ for $i=1,\ldots,n$, or just $F=z+H$.\footnote{We
should here acknowledge that in almost every other treatment of this
subject the system $F$ is written as $z-H$, which yields nicer looking
formulas for the formal inverse of $F$.  The reason for our choice is
that the formulas involving the D-Log and formal flow, which will be
developed in \S3, come out better when we write $F=z+H$.} In most
applications the power series $H=(H_1,\ldots,H_n)$ will involve will
only monomials of total degree 2 and higher, and we will often take
$H$ to be homogeneous of degree $d\ge2$.  However, these assumptions
are not necessary for what follows here.  We will associate to each
rooted tree a power series in $n$ variables based on $F$
(equivalently, on $H$).

For $T\in\bT$, a {\it labeling} of $T$ in the set $\{1,\ldots,n\}$ is
a function $f:V(T)\to\{1,\ldots,n\}$. A rooted tree $T$ with a
labeling $f$ is called a {\it labeled rooted tree,} denoted $(T,f)$.
Given such, and given $F=z+H$ as above, we make the following
definitions, for $v\in V(T)$:
\begin{enumerate}
\item $H_v=H_{f(v)}$.
\item $D_v=D_{f(v)}$.
\item $D_{v^+}=\prod_{w\in v^+} D_w$.
\item $P_{T,f}=\prod_{v\in V(T)} D_{v^+}H_v$.
\end{enumerate}
Finally, we define systems of power series
$P_T=(P_{T,1},\ldots,P_{T,n})$ and $\cP_T=(\cP_{T,1},\ldots,\cP_{T,n})$
by summing over all
labelings of $T$ having a fixed label for the root:
\begin{align}
P_{T,i}=&\sum_{\substack{{f:V(T)\to\{1,\ldots,n\}}\\f(\text{rt}_T)=i}}P_{T,f}\nno\\
\cP_{T,i}=&\,\,\frac{1}{\alpha_T}P_{T,i}\nno
\end{align}
for $i=1,\ldots,n$.

One notes that the power series $P_T$ and $\cP_T$ are dependent on the
integer $n$ and the system $H=(H_1,\ldots,H_n)\in\bC[[z]]^n$.  They
can be viewed as objects which determine functions
$\bC[[z_,\ldots,z_n]]^n\to\bC[[z_,\ldots,z_n]]^n$ for all $n\ge1$.  We
will write $P_T(H)$ and $\cP_T(H)$ when we need to emphasize this
dependence, or when we are dealing with more than one system $H$.
\vskip 10pt

\subsection{Stable Linear Independence}  We begin by establishing an
important independence property
of the objects $\{P_T\,\vert\,T\in\bT\}$.

\begin{defi}
We say rooted trees $T_1,\ldots, T_k$ are stably linear dependent if
there exist $c_1,\ldots,c_k\in\bC$  such that $\sum_{i=1}^k c_i
P_{T_i}=0$ for any integer $n\ge1$ and any homogeneous polynomial system
$H=(H_1,\ldots,H_n)$ in $n$ variables. Otherwise, we say $T_i$ are
stably linear independent.
\end{defi}

\begin{rmk}\label{homdegree}  If $H$ is homogeneous of degree $d$ and if
$T\in\bT_m$,
then $P_T(H)$ is homogeneous of degree $(d-1)m+1$.  Thus if we partition
$\{T_1,\ldots, T_k\}$ according to the number of vertices in a tree,
then $T_1,\ldots, T_k$ are stably linear independent if and only if each
partition is a stablly linearly independent set of trees.
\end{rmk}

\begin{lemma}
Suppose $\sum_{i=1}^k c_i P_{T_i}(H)=0$ for any integer $n\ge1$ and any
homogeneous polynomial system $H$ in $n$ variables. Then $\sum_{i=1}^k
c_i P_{T_i}(H)=0$ for any  system of power series $H=(H_1,\ldots,H_n)$
in $n$ variables.
\end{lemma}
\pf
We first prove it for any polynomial $H$ (not necessarily homogeneous)
in $n$ variables
by introducing a new variable $z_{n+1}$ and homogenizing $H$ using
$z_{n+1}$.  Call the resulting homogeneous system $\Bar H$. Setting
$\widetilde H=(\Bar H, H_{n+1}=0)$, we have
\BQn
\sum_{i=1}^k c_i P_{T_i}(H, z)=\sum_{i=1}^k c_i P_{T_i}(\widetilde H,
z)|_{z_{n+1}=1}=0,
\EQn
which proves the lemma for $H$ a polynomial system.  For an arbitrary
system of power series $H$ we note that if $T$ is a tree with $r$ edges,
the homogeneous summands of degree $\le d$ in $P_T(H)$ depends only on
the homogeneous summands of $H$ having degree $\le d+r$.  Taking $r$ to
be the maximum of the numbers of edges in $T_1,\ldots,T_k$, then all
terms of degree $\le d$ in $\sum_{i=1}^k c_i P_{T_i}(H)$ depend only the
homogeneous summands of $H$ having degree $\le d+r$.  Taking $\widehat
H$ to be the polynomial truncation of $H$ of degree $d+r$, we see that
$\sum_{i=1}^k c_i P_{T_i}(H)$ and $\sum_{i=1}^k c_i P_{T_i}(\widehat H)$
coincide up through degree $d$.  Since the latter is zero ($\widehat H$
being a polynomial system) and $d$ is arbitrary, we must have
$\sum_{i=1}^k c_i P_{T_i}(H)=0$.
\epfv

\begin{theo}\label{stabledepend}
Any  rooted trees $T_i$ $(i=1, 2, \cdots, k)$ with $T_i \ncong T_j$ for
any $i \neq j$ are stably linear independent.
\end{theo}

Before giving the proof we will define a polynomial system depending on
a rooted tree.  Given a rooted tree $T$ with $m$ vertices, we create
variables $z_1,\ldots,z_m$.  Label the edges $e_2,\ldots,e_m$ and assign
each variable $z_i$ with $2\le i\le m$ to the edge $e_i$.  Label the
vertices as follows: $v_1=\text{rt}_T$, and for $i=2,\ldots,m$ let $v_i$
be the vertex of $e_i$ which is furthest from the root. For each vertex
$v_i\in V(T)$, we define $H_i$ to be the product of all the variables
assigned to the edges connecting $v_i$ with its children. (Thus if $v_i$
is a leaf, we have $H_i=1$.) Set $H_T=(H_1,\ldots,H_m)$. We have:

\begin{lemma}
 Let $T$ and $T'$ be two rooted trees with same number of vertices, then
\BQn
P_{T'}(H_{T})=\begin{cases} (0,\ldots,0) \quad \quad \quad \mbox{if}
\quad T\ncong
T' \\
(\alpha_T,0,\ldots,0) \quad \mbox{if} \quad T\simeq T' \end{cases}
\EQn
\end{lemma}
\pf The following facts are not difficult to verify, and provide a
sketch of the proof:  Each coordinate $H_{T,i}$ of $H_T$ is a monomial
which is linear or constant with respect to each variable $z_i$  Each
coordinate is constant with respect to in $z_1$.  Each variable $z_i$
with $i\ge 2$ appears in precisely one coordinate $H_{T,j}$, and $i\ne
j$.  $P_{T'}(H_{T})$ is a homogeneous system of degree zero, and must be
equal to either $0$ or $1$.  If a labeling $f:V(T')\to\{1,\ldots,m\}$ is
not bijective, then $P_{T',f}=0$ since it would entail differentiating
two different coordinates $H_{T,i}, H_{T,j}$ with respect to the same
variable, or differentiating some $H_{T,i}$ twice by the same variable,
or differentiating some $H_{T,i}$ by $z_i$, all of which give zero.
Moreover, if $f(\text{rt}_{T'})\ne1$ then $P_{T',f}=0$, since it would
entail differentiation by $z_1$, and therefore $P_{T'}(H_{T})$ is zero
except possibly in the first coordinate.

With this it is not hard to show:  If $f:V(T')\to\{1,\ldots,m\}$ is a
labeling for which $P_{T',f}\ne0$, then the function $V(T')\to V(T)$
defined by $w\mapsto v_{f(w)}$ gives and isomorphism of $\varphi:T'\to
T$.  Finally, the group $\text{Aut}\,T$ acts freely and transitively on
the set of labelings $f:V(T)\to\{1,\ldots,m\}$ for which $P_{T,f}\ne0$.
The lemma follows easily from these statements.
\epfv

\noindent{\it Proof of Theorem \ref{stabledepend}}:
Suppose that $\sum_{i=1}^k c_i P_{T_i}(z)=0$ with $c_1\neq 0$. Choose
$H=H_{T_1}$, then there must exist
$j\ne 1$ such that $P_{T_j}(H_{T_1})\neq 0$. By the lemma above, we have
$T_1\simeq T_j$.
\epfv

If $H=(H_{1},\ldots,H_{n})$ is a system of power series such that
each $H_{i}$ has only terms of degree $d$ and higher, the power
series $\cP_{T}$ has only terms of degree $(d-1)v(T)+1$ and higher.
Hence if $d\ge2$ a sum of the form $\sum_{T\in\bT}c_{T}\cP_{Y}$
makes sense, since only finitely many terms contribute to any
specified homogeneous summand.  With this observation, we state the
following consequence of stable linear independence.

\begin{corol}\label{tlic} Suppose we have a collection
$\{c_{T}\}\subset\bC$ indexed by the rooted trees $T\in\bT$ such that
$\sum_{T\in\bT}c_{T}\cP_{T}=0$ for for any integer $n\ge1$ and any
system of power series $H=(H_{1},\ldots,H_{n})$ with $H$ having only
terms of degree $\ge2$.  Then $c_{T}=0$ for all $T\in\bT$.
\end{corol}
\pf We consider systems $H$ which are homogeneous polynomial systems
of degree $d\ge2$.  In this case $\cP_{T}$ is homogeneous of degree
$(d-1)v(T)+1$, so the homogeneous summands of
$\sum_{T\in\bT}c_{T}\cP_{T}$ are the finite sums
$\sum_{T\in\bT_{N}}c_{T}\cP_{T}$ for $N\in\bN$, so these must be zero.
By Theorem \ref{stabledepend} applied to the finite set of trees
$\bT_{N}$, we must have $c_{T}=0$ for all $T\in\bT_{N}$.  \epfv

Recall that we are writing $D_i$ for the operator
$\frac{\partial}{\partial z_i}$.  We will denote by $D$ the column
vector $(D_1,\ldots,D_n)^{\text{t}}$.  We now define a differential
operator on $\bC[[z]]$ for each $T\in\bT$.

\begin{defi}\label{Ddef} For $T\in\bT$, we denote by $D_T$ the
differential operator $P_T D=\sum_{i=1}^n P_{T,i}D_i$.  We will write
$\cD_T$ for the operator $\cP_T D=\frac{1}{\alpha_T}D_T$.
\end{defi} \vskip 10pt

\subsection{Tree Surgery} We will now discuss some ``surgical"
procedures on trees.  Given $T\in\bT$ and $e\in E(T)$, the removal of
the edge $e$ from $T$
gives a disconnected graph with two connected components which are
trees.  We denote by $T_e$ the component containing $\text{rt}_T$, and
by $T'_e$ the other component.  We give $T_e$ and $T'_e$ the structure
of rooted trees by setting $\text{rt}_{T_e}=\text{rt}_T$ and
$\text{rt}_{T'_e}=v'_e$.

Given rooted trees $T$ and $T'$ and $v\in V(T)$, we denote by
$$T'\multimap_v T$$ the tree obtained by connecting $\text rt_{T'}$
and $v$ by a newly created edge, and setting
$\text{rt}_{(T'\multimap_v T)}=\text{rt}_T$.  We will refer to the
newly created edge as the {\it connection edge} of $T'\multimap_v T$.
Note that for and tree $T$ and edge $e\in E(T)$ we have an obvious
isomorphism $T\simeq(T'_e\multimap_{v_e} T_e)$ which is the identity on
$T_e$ and $T'_e$.

Given $e,f\in E(T)$, we say ``$f$ lies below $e$", and write $e\succ
f$, if $f\in E(T_e)$.  This merely says that $f$ remains when we
``strip away" $e$ and $T'_e$.  One can easily see that this relation
is not transitive.  However, if we write $$e_1\succ \cdots \succ
e_r\,,$$ for $e_1,\ldots,e_r\in E(T)$, we will mean by this that
$e_i\succ e_j$ if $i<j$.

A sequence $\vec{e}=(e_1,\ldots,e_r)\in E(T)^r$ with $e_1\succ \cdots
\succ e_r$ determines a sequence of subtrees $T_{\vec{e},1},\ldots
T_{\vec{e},r+1}$ as follows: Setting $T_{\vec{e},1}=T'_{e_1}$ and let
$S_2=T_{e_1}$, noting that $e_2,\ldots,e_r\in E(S_2)$.  For
$i=1,\ldots,r$, assume $T_{\vec{e},1},\ldots,T_{\vec{e},i-1},S_i$ are
defined with $e_i,\ldots,e_r\in E(S_i)$.  Set
$T_{\vec{e},i}=(S_{i})'_{e_i}$ and $S_{i+1}=(S_{i})_{e_i}$.  Finally,
set $T_{\vec{e},r+1}=S_{r+1}$.

For any integer $r\geq 1$ and $T\in\bT$, create an indeterminate
$Y_T^{(r)}$.  Denote this set of variables (for all $T$ and $r$) by
$Y$.  Extend the action of the operators $D_T$ and $\cD_T$ to
$\bC[[z]][Y]$ by making each indeterminate of $Y$ a constant.

\begin{lemma}\label{key-lemma}  Let $r,m\ge1$ be integers and
$S\in\bT$.  Then
\allowdisplaybreaks{
\BQ
&{}&
\sum_{\substack{(T_1,\ldots,T_r)\in\bT^r\\v(T_1)+\cdots
+v(T_r)+v(S)=m}}\left[Y_{T_1}^{(1)}\cD_{T_1}\right]
\cdots \left[Y_{T_r}^{(r)}\cD_{T_{r}}\right]\cP_S \nno\\
&{}&\qquad\qquad\qquad\qquad =
\sum_{T\in {\bT}_m}\,\,\,\,
\sum_{\substack{\vec{e}=(e_1,\ldots,e_r)\in E(T)^r\\e_1\succ \cdots
\succ e_r\\T_{\vec{e},r+1}\simeq S}}Y_{T_{\vec{e},1}}^{(1)}\cdots
Y_{T_{\vec{e},r}}^{(r)}\cP_T \label{key-eq}
\EQ}
\end{lemma}
\pf  Note that both sums are finite, so the expression makes sense for
any $H\in\bC[[z]]^n$.

We first consider the case $r=1$. For $T'\in\bT$ we have
$$D_{T'}P_S=\sum_{v\in V(S)}P_{(T'\multimap_v S)}\,.$$ Hence
\allowdisplaybreaks{
\begin{align}
\sum_{\substack{T'\in\bT\\v(T')+v(S)=m}}Y_{T'}^{(1)}D_{T'}P_S
&=\sum_{\substack{T'\in\bT\\v(T')+v(S)=m}}\sum_{v\in
V(S)}Y_{T'}^{(1)}P_{(T'\multimap_v S)}\nno\\
&=\sum_{T\in\bT_m}\sum_{T'\in\bT}\sum_{\substack{v\in
V(S)\\(T'\multimap_v S)\simeq T}}Y_{T'}^{(1)}P_{(T'\multimap_v S)}
\,.\nno\end{align}
For a fixed $T\in\bT_m$ we wish to count the occurrences of $P_T$ in the
last expression.  Toward this end, for $T'\in\bT$ let
\begin{align}
I_{T,T',S}&=\{v\in V(S)\,\vert\,(T'\multimap_v S)\simeq T\}\nno\\
J_{T,T',S}&=\{\bar e\in E(T)/\text{Aut}\,(T)\,\vert\,T'_e\simeq
T',T_e\simeq S\text{ (for any $e$ representing $\bar e$)}\}\nno
\end{align}}
We will define as function $\Phi:I_{T,T',S}\to J_{T,T',S}$ as follows:
Given $v\in I_{T,T',S}$, choose an isomorphism $\varphi:(T'\multimap_v
S)\overset{\simeq}{\longrightarrow}T$, and let $e$ be the image under
$\varphi$ of the connection edge in $T'\multimap_v S$.  Letting $\bar e$
be the class of $e$ in $E(T)/\text{Aut}\,(T)$, we clearly have $\bar
e\in J_{T,T',S}$.  To see that $\bar e$ is independent of the choice of
$\varphi$, suppose $\gamma:(T'\multimap_v
S)\overset{\simeq}{\longrightarrow}T$ sends the connection edge to $f\in
E(T)$.  Then $\gamma\varphi^{-1}(e)=f$, hence $\bar f=\bar e$ in
$E(T)/\text{Aut}\,(T)$.  Therefore we have a well-defined function
$\Phi$, which is obviously surjective.

We claim that for $v\in I_{T,T',S}$ the orbit of $v$ under
$\text{Aut}\,(S)$ is precisely the fiber of $v$ under $\Phi$.  It is
clear that if $w\sim v$ by the action of $\text{Aut}\,(S)$ then
$(T'\multimap_w S)\simeq (T'\multimap_v S)\simeq T$, with the first
isomorphism taking one connection edge to the other, which shows $w\in
I_{T,T',S}$.  Choosing appropriate isomorphisms $(T'\multimap_w
S)\overset{\rho}{\longrightarrow} (T'\multimap_v
S)\overset{\varphi}{\longrightarrow}T$, we see that the image $e$ of the
connection edge of $T'\multimap_v S$ under $\varphi$ is also the image
of the connection edge of $T'\multimap_w S$ under $\varphi\rho$, hence
$\Phi(w)=\Phi(v)$.  Moreover, if $w\in I_{T,T',S}$ is \underbar{any}
element for which $\Phi(w)=\Phi(v)$, then we have isomorphisms
$$(T'\multimap_w S)\overset{\gamma}{\longrightarrow} T
\overset{\varphi}{\longleftarrow}(T'\multimap_v S)$$ such that the same
$e\in E(T)$ is the image of both connection edges.  (This can be
achieved after modifying by an automorphism of $T$.)  It follows that
$\gamma^{-1}\varphi:(T'\multimap_v S)\to
(T'\multimap_w S)$ carries one connection edge to the other, so it
restricts to an automorphism of $S$ sending $v$ to $w$.  Hence $w\sim
v$.  Therefore the above sum can be written as
\allowdisplaybreaks{
\begin{align}
&\sum_{T\in\bT_m}\,\,\sum_{T'\in\bT}\,\,\sum_{v\in
I_{T,T',S}}Y_{T'}^{(1)}P_T\nno\\
&\quad\quad=\sum_{T\in\bT_m}\,\,\sum_{T'\in\bT}\,\,\sum_{\bar e\in
J_{T,T',S}}s_{T_e}(v_e)\,Y_{T'}^{(1)}P_T\nno\\
&\quad\quad=\sum_{T\in\bT_m}\,\,\sum_{\substack{\bar e\in
E(T)/\text{Aut}\,(T)\\T_e\simeq S}}s_{T_e}(v_e)\,Y_{T'_e}^{(1)}P_T\nno
\end{align}}
where $s_{T_e}(v_e)$ is the orbit size of $v_e$ under the action of
$\text{Aut}\,T_e$, for some (any) $e\in E(T)$ representing $\bar e$. The
number of edges representing $\bar e$ is $\alpha_T /\alpha_{T,e}$, hence
the inner sum can be altered to run over all $e\in E(T)$ at the cost of
dividing by $\alpha_T /\alpha_{T,e}$, yielding
$$\sum_{T\in\bT_m}\frac{1}{\alpha_T}\sum_{\substack{e\in E(T)\\T_e\simeq
S}}\alpha_{T,e}\,s_{T_e}(v_e)\,Y_{T'_e}^{(1)}P_T\,.$$
An automorphism of $T$ fixing $e\in E(T)$ restricts to an automorphism
of $T'_e$ and an automorphism of $T_e$ fixing $v_e$.  Conversely, given
the latter pair we get a unique automorphism of of $T$ preserving $e$.
It follows that $\alpha_{T,e}=\alpha_{T'_e}\alpha_{T_e,v_e}$.  Also we
have $s_{T_e,v_e}=\alpha_{T_e}/\alpha_{T_e,v_e}$. Incorporating these
facts and putting together the above equalities, we get
$$\sum_{\substack{T'\in\bT\\v(T')+v(S)=m}}Y_{T'}^{(1)}D_{T'}P_S
=\sum_{T\in\bT_m}\frac{1}{\alpha_T}\sum_{\substack{e\in
E(T)\\T_e\simeq S}}\alpha_{T'_e}\alpha_{S}Y_{T'_e}^{(1)}P_T$$
Dividing the equation by $\alpha_{S}$ and substituting
$\frac{1}{\alpha_R}Y_R^{(1)}$ for $Y_R^{(1)}$ for each $R\in\bT$ yields
$$\sum_{\substack{T'\in\bT\\v(T')+v(S)=m}}Y_{T'}^{(1)}\cD_{T'}\cP_S
=\sum_{T\in\bT_m}\,\sum_{\substack{e\in
E(T)\\T_e\simeq S}}Y_{T'_e}^{(1)}\cP_T\,,$$
which is precisely the assertion of the lemma for $r=1$.

For $r\ge2$ we apply induction as follows:
\allowdisplaybreaks{
\begin{align}
&\sum_{\substack{(T_1,\ldots,T_r)\in\bT^r\\v(T_1)+\cdots
+v(T_r)+v(S)=m}}\left[Y_{T_1}^{(1)}\cD_{T_1}\right]
\cdots \left[Y_{T_r}^{(r)}\cD_{T_{r}}\right]\cP_S\nno\\
&\quad
=\sum_{T_1\in\bT}Y_{T_1}^{(1)}\cD_{T_1}\sum_{\substack{(T_2,\ldots,T_r)
\in\bT^{r-1}\\v(T_2)+\cdots
+v(T_r)+v(S)=m-v(T_1)}}\left[Y_{T_2}^{(2)}\cD_{T_2}\right]
\cdots \left[Y_{T_r}^{(r)}\cD_{T_{r}}\right]\cP_S\,.\nno
\end{align}}
Applying induction and a substitution of variables $Y_t^{(i+1)}$ for
$Y_t^{(i)}$ to the inner sum, this equals
{\allowdisplaybreaks\begin{align}
&\sum_{T_1\in\bT}Y_{T_1}^{(1)}\cD_{T_1}\sum_{R\in
{\bT}_{m-v(T_1)}}\,\,\,\,
\sum_{\substack{\vec{e}=(e_1,\ldots,e_{r-1})\in E(R)^{r-1}\\e_1\succ
\cdots
\succ e_{r-1}\\T_{\vec{e},r}\simeq S}}Y_{T_{\vec{e},1}}^{(2)}\cdots
Y_{T_{\vec{e},r-1}}^{(r)}\cP_R
\nno\\
=&\sum_{\substack{T_1,R\in\bT\\v(T_1)+v(R)=m}}Y_{T_1}^{(1)}\cD_{T_1}\,\,\,\,
\sum_{\substack{\vec{e}=(e_1,\ldots,e_{r-1})\in
E(R)^{r-1}\\e_1\succ \cdots \succ
e_{r-1}\\T_{\vec{e},r}\simeq S}}Y_{T_{\vec{e},1}}^{(2)}\cdots
Y_{T_{\vec{e},r-1}}^{(r)}\cP_R\nno\\
=&\sum_{R\in\bT}\left[\sum_{\substack{T_1\in\bT\\v(T_1)+v(R)=m}}Y_{T_1}^{(1)}\cD_{T_1}\cP_R\right]\,\,\,\,\sum_{\substack{\vec{e}=(e_1,\ldots,e_{r-1})\in
E(R)^{r-1}\\e_1\succ \cdots
\succ e_{r-1}\\T_{\vec{e},r}\simeq S}}Y_{T_{\vec{e},1}}^{(2)}\cdots
Y_{T_{\vec{e},r-1}}^{(r)}\,.\nno
\end{align}}
Now we apply the case $r=1$ to the bracketed expression to obtain
{\allowdisplaybreaks
\begin{align}
&\sum_{R\in\bT}\left[\sum_{T\in\bT_m}\,\sum_{\substack{e\in
E(T)\\T_e\simeq
R}}Y_{T'_e}^{(1)}\cP_T\right]\,\,\,\,\sum_{\substack{\vec{e}=(e_1,\ldots,e_{r-1})\in
E(R)^{r-1}\\e_1\succ \cdots
\succ e_{r-1}\\T_{\vec{e},r}\simeq S}}Y_{T_{\vec{e},1}}^{(2)}\cdots
Y_{T_{\vec{e},r-1}}^{(r)}\nno\\
&\sum_{T\in\bT_m}\sum_{R\in\bT}\sum_{\substack{e\in
E(T)\\T_e\simeq
R}}\,\,\,\,\sum_{\substack{\vec{e}=(e_1,\ldots,e_{r-1})\in
E(R)^{r-1}\\e_1\succ \cdots
\succ
e_{r-1}\\T_{\vec{e},r}\simeq
S}}Y_{T'_e}^{(1)}Y_{T_{\vec{e},1}}^{(2)}\cdots
Y_{T_{\vec{e},r-1}}^{(r)}\cP_T\nno\\
&=\sum_{T\in {\bT}_m}\,\,\,\,
\sum_{\substack{\vec{e}=(e_1,\ldots,e_r)\in E(T)^r\\e_1\succ \cdots
\succ e_r\\T_{\vec{e},r+1}\simeq S}}Y_{T_{\vec{e},1}}^{(1)}\cdots
Y_{T_{\vec{e},r}}^{(r)}\cP_T\nno
\end{align}
}
which completes the proof.
\epfv

Suppose the system of power series $H=(H_{1},\ldots,H_{n})$ has the
property that each $H_{i}$ involves only monomials of degree $\ge2$ in
$z_{1},\ldots,z_{n}$.  Then one easily verifies that for $T\in\bT$,
$\cP_{T}$ involves only monomials of degree $\ge v(T) +1$.  It follows
that for a monomial $M$ in $z$ of degree $m$, $\cD_{T}\cdot M$
involves only monomials of degree $\ge m+v(T)$.  Therefore infinite
sums such as $\sum_{T\in\bT}\cP_{T}$ and $\sum_{T\in\bT}\cD_{T}$ make
sense in this situation.  The following two corollaries of Lemma
\ref{key-lemma} are based on this observation.  The equations in both
corollaries take place in the ring $\bC[Y][[z]]$, where $Y$ represents
the infinite set of variables $\{Y_{T}^{(i)}\,|\,T\in\bT,i\in\bZ^{+}\}$.

\begin{corol}\label{key-cor1}  Suppose the system of power series $H$
involves only monomials of degree $\ge2$.
Let $r\ge1$ be an integer and
$S\in\bT$.  Then
\allowdisplaybreaks{
\BQ
&{}&
\sum_{(T_1,\ldots,T_r)\in\bT^r}\left[Y_{T_1}^{(1)}\cD_{T_1}\right]
\cdots \left[Y_{T_r}^{(r)}\cD_{T_{r}}\right]\cP_S \nno\\
&{}&\qquad\qquad\qquad =
\sum_{T\in\bT}\,\,\,\,
\sum_{\substack{\vec{e}=(e_1,\ldots,e_r)\in E(T)^r\\e_1\succ \cdots
\succ e_r\\T_{\vec{e},r+1}=S}}Y_{T_{\vec{e},1}}^{(1)}\cdots
Y_{T_{\vec{e},r}}^{(r)}\cP_T \label{key-eq1}
\EQ}
\end{corol}
\pf  We simply sum (\ref{key-eq}) over all $m\ge1$, noting the
convergence of the sums by the observations above.
\epfv

\begin{corol}\label{key-cor}  Suppose the system of power series $H$
involves only monomials of degree $\ge2$.  Let $k\ge2$ be an integer.
Then
\BQ
&{}&
\sum_{(T_1,\ldots,T_k)\in\bT^k}\left[Y_{T_1}^{(1)}\cD_{T_1}\right]
\cdots
\left[Y_{T_{k-1}}^{(k-1)}\cD_{T_{k-1}}\right]\left[Y_{T_k}^{(k)}\cP_{T_{k}}\right]
\nno\\
&{}&\qquad\qquad =
\sum_{\substack{{T\in\bT}\\v(T)\ge2}}\,\,\,\,
\sum_{\substack{\vec{e}=(e_1,\ldots,e_{k-1})\in E(T)^{k-1}\\e_1\succ
\cdots
\succ e_{k-1}}}Y_{T_{\vec{e},1}}^{(1)}\cdots
Y_{T_{\vec{e},k}}^{(k)}\cP_T \label{key-eq2}
\EQ
\end{corol}
\pf We apply Corollary \ref{key-cor1}, multiplying both sides of
(\ref{key-eq1}) by $Y_{S}^{(r+1)}$, setting $k=r+1$, summing over all
$S\in\bT$.  Note that the singleton tree contributes 0 in
(\ref{key-eq1}) for any $r\ge1$, and thus the qualifier $v(T)\ge2$ in
(\ref{key-eq2}).  \epfv

\renewcommand{\theequation}{\thesection.\arabic{equation}}
\renewcommand{\therema}{\thesection.\arabic{rema}}
\setcounter{equation}{0}
\setcounter{rema}{0}

\section{D-Log and Formal Flow}

We will henceforth be restricting our attention to systems of power
series $F=(F_1, \ldots,F_n)\in\bC[[z]]^n$ of the form $F_i=z_i+H_i$
with $H_i$ involving only monomials of degree 2 and higher, for
$i=1,\ldots,n$.  We refer to this condition by saying ``$F$ is of the
form {\it identity plus higher.}'' Such a system determines a
$\bC$-algebra automorphism of $\bC[[z]]$, namely the automorphism
which sends $z_i$ to $F_i$ for $i=1,\ldots,n$.

\subsection{The D-Log} The following proposition appears as
Proposition 2.1 in \cite{Z}.
\begin{propo}\label{Prpo-3.1}
For any $F=(F_1, F_2, \cdots , F_n)\in\bC[[z]]^n$ of the form identity
plus higher, there exists a
unique system of power series
$$a=(a_1, a_2, \cdots , a_n)\in{\bC}[[z]]^n$$ involving only monomials
of
degree 2 and higher such that, letting $A=aD=\sum_{i=1}^n a_iD_i$, we
have
\BQ\label{EF1}
\exp(A)\cdot z = F
\EQ
where
$$
\exp(A)=\sum_{k=0}^\infty \frac {A^k}{k!}
$$
and $z=(z_1,\ldots,z_n)$.
\end{propo}

The reader will easily verify that the infinite sum $\exp(A)\cdot Q$
makes sense for any $Q\in\bC[[z]]^n$ due to the fact that, for any
integer $d\ge0$, only finitely many terms $\frac {A^k}{k!}\cdot Q$
contribute to the degree $d$ homogeneous summand.  This is due to the
fact that $a$ involves only terms of degree 2 and higher.

\begin{rmk} \label{deraut} It is well known that the exponential of a
derivation on any $\bQ$-algebra, when it makes sense, is an
automorphism of that algebra.  Any subring lying in the kernel of the
derivation will be fixed by this automorphism.  It follows from this
fact, the comment
above, and Proposition
\ref{Prpo-3.1} that $\exp(A)$ is the $\bC$-algebra automorphism of
$\bC[[z]]$
which sends $z_i$ to $F_{i}$, for $i=1,\ldots,n$. \end{rmk}

\begin{defi}
We call the unique system of power series $a=(a_1,\dots,a_n)$ obtained
above the {\it Differential Log} or {\it
D-Log} of the formal system $F$.
\end{defi}

\subsection{Coefficients $\phi_T$ of the D-Log}

\begin{theo}\label{Theo-3.3}
There exists a unique set of rational numbers $\{\phi_T\}$ indexed by
the set of rooted trees $T\in \bT$,
such that
\BQ\label{E-3.4}
a=\sum_{T\in \bT} \phi_T {\mathcal P_T}\,.
\EQ
These numbers satisfy, and are uniquely determined by, the following
properties:
\begin{align}
\phi_T&=1 \text{ when $v(T)=1$ (i.e., $T=\circ$, the
singleton tree)}\notag\\
\phi_T&=-\sum_{k=2}^{v(T)} \frac 1{k!}
\sum_{\substack{\vec e=(e_1,\ldots,e_{k-1})\in E(T)^{k-1}\\e_1\succ
\cdots\succ e_{k-1}}}
\phi_{T_{\vec e,1}}\phi_{T_{\vec e,2}}\cdots\phi_{T_{\vec e,k}} \text{
when $v(T)\ge2$}\,.\label{RecurForPhi}
\end{align}
The latter formula can be restated as:
\BQ \label{R00}
\sum_{k=1}^{v(T)} \frac 1{k!}
\sum_{\substack{\vec e=(e_1,\ldots,e_{k-1})\in E(T)^{k-1}\\e_1\succ
\cdots\succ e_{k-1}}}
\phi_{T_{\vec e,1}}\phi_{T_{\vec e,2}}\cdots\phi_{T_{\vec e,k}}=0\,.
\EQ
(Here we must interpret the $k=1$ summand as $\phi_T$.)
\end{theo}
\pf Let us define $\phi_T$ by (\ref{RecurForPhi}) and set
$a'=\sum_{T\in\bT}\phi_T\cP_T$, $A'=a'D$. Then $A'=
\sum_{T\in\bT}\phi_T\cD_T$.  We have
\allowdisplaybreaks{
\begin{align}
\exp(A')\cdot z&=\sum_{k=0}^\infty \frac{{A'}^k}{k!}\cdot z\nno\\
&=\sum_{k=0}^\infty
\frac{1}{k!}\left(\sum_{T\in\bT}\phi_T\cD_{T}\right)^{k}\cdot z\nno\\
&=z+\sum_{k=1}^\infty
\frac{1}{k!}\sum_{(T_1,\ldots,T_k)\in\bT^k}\left[\phi_{T_1}\cD_{T_1}\right]\cdots
\left[\phi_{T_k}\cD_{T_k}\right]\cdot z\nno\\
&=z+\sum_{T\in\bT}\phi_T\cD_T\cdot z + \sum_{k=2}^\infty
\frac{1}{k!}\sum_{(T_1,\ldots,T_k)\in\bT^k}\left[\phi_{T_1}\cD_{T_1}\right]\cdots
\left[\phi_{T_k}\cD_{T_k}\right]\cdot z\nno\\
\intertext{Using the fact that $\cD_T\cdot z=\cP_T$:}
=z+&\sum_{T\in\bT}\phi_T\cP
+ \sum_{k=2}^\infty
\frac{1}{k!}\sum_{(T_1,\ldots,T_k)\in\bT^k}\left[\phi_{T_1}\cD_{T_1}\right]\cdots
\left[\phi_{T_{k-1}}\cD_{T_{k-1}}\right]\left[\phi_{T_k}\cP_{T_k}\right]\nno\\
\intertext{Applying Corollary \ref{key-cor}, substituting
$Y_{T_i}^{(i)}=\phi_{T_i}$:}
=z+&\sum_{T\in\bT}\phi_T\cP_T
+\sum_{k=2}^\infty
\frac{1}{k!}\sum_{\substack{T\in\bT\\v(T)\ge2}}\,\,\,\sum_{\substack{\vec{e}=(e_1,\ldots,e_{k-1})\in
E(T)^{k-1}\\e_1\succ \cdots \succ e_{k-1}}}\phi_{T_{\vec
e,1}}\cdots\phi_{T_{\vec e,k}}\cP_T\nno\\
\intertext{Letting $S$ be the singleton tree:}
=z+&\phi_S\cP_S
+\sum_{\substack{T\in\bT\\v(T)\ge2}}\left(\sum_{k=1}^{v(T)}\frac{1}{k!}
\sum_{\substack{\vec{e}=(e_1,\ldots,e_{k-1})\in E(T)^{k-1}\\e_1\succ
\cdots \succ e_{k-1}}}\phi_{T_{\vec e,1}}\cdots\phi_{T_{\vec
e,k}}\right)\cP_T\nno\\
\intertext{Since $\cP_S=H$, and, by definition, $\phi_S=1$ and the sum
in parentheses is 0:}
&=z+H=F\nno
\end{align}
By the uniqueness property of $a$ we must have $a'=a$.  The uniqueness
of the expression (\ref{E-3.4}) for $a$ follows from Theorem
\ref{stabledepend} \epfv}

\noindent {\bf Chains and Shrubs.} Two special types of trees are the
``chains" and the ``shrubs", mentioned in \S2.  Given and integer
$n\ge1$ we let
$C_n\in\bT_n$ be the {\it chain} with $n$ vertices, which is the unique
rooted tree in $\bT_n$ of height $n-1$.  For $n\ge 0$ we let
$S_n\in\bT_{n+1}$ be the {\it shrub} with $n+1$ vertices, which is the
unique rooted tree in $\bT_{n+1}$ of height $\le1$ (Equality holds
unless $n=0$.).  Note that $C_1=S_0=\circ\,$, the singleton tree.

By using the recurrence formula (\ref{RecurForPhi}), we can calculate
$\phi_T$ for chains and shrubs as follows.  Consider the generating
functions
\begin{align}
c(x)&=\sum_{n=1}^\infty \phi_{C_n}x^n\nno\\
s(x)&=\sum_{n=0}^\infty \phi_{S_n}\frac{x^n}{n!}\nno
\end{align}
Then we have:

\begin{corol}\label{corol-3.4}
The generating functions $c(x)$ and $s(x)$ are given by:
\begin{align}
\text{\rm (a)}\qquad c(x)&=\ln (1+x)\label{chainfct}\\
\text{\rm (b)}\qquad s(x)&=\frac{x}{e^{x}-1}\label{shrubfct}
\end{align}
In particular, we have $\phi_{C_n}=(-1)^{n-1}\frac{1}n$ for all $n\ge1$
and
$\phi_{S_n}=b_n$, where $b_0,b_1,b_2,\ldots\,$ are
the Bernoulli numbers$\,$\footnote{This indexing and signage differs
from an alternate definition of the Bernoulli numbers as the
sequence $B_1,B_2,\ldots\,$ defined by
$$\frac{x}{e^{x}-1}=1-\frac{1}{2}x+\sum_{n=1}^\infty
(-1)^{n-1}\frac{B_n}{(2n)!}x^{2n}\,.$$ Thus the relationship is
$B_n=(-1)^{n-1}b_{2n}$ for $n\ge1$.} defined by
$\frac{x}{e^{x}-1}=\sum_{n=0}^\infty b_n \frac{x^n}{n!}$.
\end{corol}
\pf
(a) According to (\ref{RecurForPhi}) we have
{\allowdisplaybreaks\begin{align}
c(x)&=x-\sum_{n=2}^\infty
\left(\sum_{k=2}^{v(C_n)}  \frac 1{k!}
\sum_{\substack{\vec e=(e_1,\ldots,e_{k-1})\in E(C_n)^{k-1}\\e_1\succ
\cdots\succ e_{k-1}}}
\phi_{T_{\vec e,1}}\phi_{T_{\vec e,2}}\cdots\phi_{T_{\vec e,k}}\right )
x^n\nno\\
\intertext{Noting that $v(C_n)=n$ and each $T_{\vec e,j}$ is also a
path:}
&=x-\sum_{n=2}^\infty
\sum_{k=2}^{n} \frac 1{k!}
\sum_{\substack{(i_1,\ldots,i_k)\in\bN^k\\i_1+\cdots+i_k=n}}\prod_{j=1}^{k}
\phi_{C_{i_j}}x^{i_j}\nno\\
&=x-\sum_{n=2}^\infty
\left(\sum_{k=2}^{n} \frac 1{k!}
(\text{coefficient of $x^n$ in $c(x)^k$})\right)\,x^n\nno\\
&=x-\sum_{n=2}^\infty
\left(\text{coefficient of $x^n$ in $\sum_{k=2}^{n} \frac
1{k!}c(x)^k$}\right)\,x^n\nno\\
&=x-\sum_{n=2}^\infty
\left(\text{coefficient of $x^n$ in $\sum_{k=2}^{\infty} \frac
1{k!}c(x)^k$}\right)\,x^n\nno\\
&=x-\sum_{n=2}^\infty\frac 1{n!}c(x)^n\nno\\
&=x-(e^{c(x)}-c(x)-1)\nno
\end{align}}
Solving for $c(x)$ in the equation $c(x)=x-(e^{c(x)}-c(x)-1)$ gives
(\ref{chainfct}). \vskip 10pt

(b) Again by (\ref{RecurForPhi}) we have
{\allowdisplaybreaks\begin{align}
s(x)&=1-\sum_{n=1}^\infty
\left(\sum_{k=2}^{v(S_n)}\frac 1{k}
\sum_{\substack{\vec e=(e_1,\ldots,e_{k-1})\in E(S_n)^{k-1}\\e_1\succ
\cdots\succ e_{k-1}}}
\phi_{T_{\vec e,1}}\phi_{T_{\vec e,2}}\cdots\phi_{T_{\vec e,k}}\right )
\frac{x^n}{n!}\nno\\
\intertext{Noting that $v(S_n)=n+1$ and precisely one $T_{\vec e,j}$ is
a shrub with all others being singletons:}
&=1-\sum_{n=1}^\infty
\left(\sum_{k=2}^{n+1} \frac 1{k!}(k-1)!\binom{n}{k-1}
\phi_{S_{n-k+1}}\right)\frac{x^n}{n!}\nno\\
&=1-x^{-1}\sum_{n=1}^\infty
\sum_{k=2}^{n+1}
\phi_{S_{n-k+1}} \frac{x^{n-k+1}}{(n-k+1)!}  \frac {x^{k}}{k!}\\
&=1-x^{-1}\left(\sum_{r=0}^\infty\phi_{S_r}\frac{x^r}{r!}\right)
\left(\sum_{s=2}^\infty\frac{x^s}{s!}\right)\nno\\
&=1-x^{-1}s(x) (e^x-x-1)
\end{align}}
Solving for $s(x)$ in the equation $s(x)=1-x^{-1}s(x) (e^x-x-1)$ gives
(\ref{shrubfct}).
\epfv

\subsection{Polynomial Coefficients $\psi_T(t)$ of Formal Flow}

Let us first recall the formal flow $F_t=\exp(tA)\cdot z$
and some of its properties. See \cite{Z} for more details.

\begin{defi}
Given an indeterminate
$t$, define the system $F_t\in\bC[t][[z]]^n$ by
\BQ
F_t=\exp(tA)\cdot z\,.
\EQ
It is called the {\rm formal flow} generated by $F$.
\end{defi}

It is easy to verify that $F_t\in\bC[t][[z]]^n$.  Therefore a
specialization $t=\alpha$, for any $\alpha\in\bC$ (or $\alpha$ in any
$\bC$-algebra), makes sense.  According to Proposition
(\ref{Prpo-3.1}), setting $t=1$ in $F_t$ recovers $F$.

The following proposition shows that $t$ behaves like an exponent
for $F$.

\begin{propo}\label{exp}  Let $t$ and $s$ be indeterminates.  Then
    \BQ
    F_{s+t}=F_{t}\circ F_{s}\nno
    \EQ
Hence setting $t=n$ in $F_{t}$, for $n\in\bN$, gives the $n$-fold
composition $F\circ\cdots\circ F$, and setting $t=-n$ gives the
$n$-fold composition $F^{-1}\circ\cdots\circ F^{-1}$ of the formal
inverse.  In particular,
$$F_{t}\vert_{{}_{t=-1}}=F^{-1}\,.$$
\end{propo}
\pf
We have
\begin{align}
    F_{s+t}&=\exp((s+t)A)\cdot z=\exp(sA+tA)\cdot z\nno\\
    &=\exp(sA)\cdot\exp(tA)\cdot z=\exp(sA)\cdot F_{t}\nno\\
\intertext{We use that fact that $exp(sA)$ is a $\bC$-algebra automorphism of
$\bC[s,t][[z]]$ (see Remark \ref{deraut}):}
    &=F_{t}(\exp(sA\cdot z))=F_{t}(F_{s})\nno\\
    &=F_{t}\circ F_{s}\nno
\end{align} \epfv

Thus $F_{t}$ can be viewed as the ``formal $t^{\text{th}}$ 
power of $F$''.

The system $F_{t}$ can be expressed in terms of the tree
expressions $\cP_{T}$ as follows:

\begin{theo}\label{Theo-3.6}
There exists a unique set of polynomials $\{\psi_T(t)\in\bQ[t]\}$
indexed by the set of rooted trees $T\in\bT$
such that
\BQ\label{compare}
F_{t}=z+
\sum_{T\in\bT} \psi_T(t){\mathcal P}_T\,.
\EQ
These polynomials are given by the formula
\BQ\label{DefForPsi}
\psi_T(t)=\sum_{k=1}^{v(T)} \frac {t^{k}}{k!}
\sum_{\substack{\vec e=(e_1,\ldots,e_{k-1})\in E(T)^{k-1}\\e_1\succ
\cdots\succ e_{k-1}}}
\phi_{T_{\vec e,1}}\phi_{T_{\vec e,2}}\cdots\phi_{T_{\vec e,k}}
\EQ
(Again we must interpret the $k=1$ summand as $\phi_T$.)
\end{theo}
\pf According to Theorem \ref{Theo-3.3} the D-Log of $F$ is given by
$a=\sum_{T\in\bT}\phi_{t}\cP_{T}$, hence we have
$A=aD=\sum_{T\in\bT}\phi_{t}\cP_{T}D=\sum_{T\in\bT}\phi_{T}\cD_{T}$
(see Definition \ref{Ddef}).  Therefore
\allowdisplaybreaks{
\begin{align}
   F_t&=\exp(tA)\cdot z=\sum_{k=0}^{\infty}\frac{t^{k}}{k!}A^{k}\cdot
z\nno\\
&=z+\sum_{k=1}^{\infty}\frac{t^{k}}{k!}\left(\sum_{T\in\bT}\phi_{T}\cD_{T}
   \right)^{k}\cdot z\nno\\
&=z+\sum_{k=1}^{\infty}\frac{t^{k}}{k!}\sum_{(T_{1},\ldots,T_{k})\in\bT}
   [\phi_{T_1}\cD_{T_1}]\cdots[\phi_{T_k}\cD_{T_k}]\cdot z\nno\\
&=z+\sum_{k=1}^{\infty}\frac{t^{k}}{k!}\sum_{(T_{1},\ldots,T_{k})\in\bT}
[\phi_{T_1}\cD_{T_1}]\cdots[\phi_{T_{k-1}}\cD_{T_{k-1}}][\phi_{T_k}\cP_{T_k}]\nno\\
   \intertext{Now we apply Corollary \ref{key-cor} to the $k\ge2$
   summands:}
   &=z+\sum_{k=1}^{\infty}\frac{t^{k}}{k!}\sum_{T\in\bT}\,\,\,
   \sum_{\substack{\vec{e}=(e_1,\ldots,e_{k-1})\in E(T)^{k-1}\\e_1\succ
\cdots
   \succ e_{k-1}}}\phi_{T_{\vec{e},1}}\cdots
   \phi_{T_{\vec{e},k}}\cP_T \nno\\
   &=z+\sum_{T\in\bT}\,\left(\sum_{k=1}^{v(T)}\frac{t^{k}}{k!}\,\,\,
   \sum_{\substack{\vec{e}=(e_1,\ldots,e_{k-1})\in E(T)^{k-1}\\e_1\succ
\cdots
   \succ e_{k-1}}}\phi_{T_{\vec{e},1}}\cdots
   \phi_{T_{\vec{e},k}}\right)\cP_T \nno
   \end{align}
This gives the desired result.  The uniqueness of $\psi_{T}$ follows
from applying stable linear independence (Corollary \ref{tlic}) to each
power of $t$ in (\ref{compare}).  \epfv}

\begin{lemma}\label{lemma-3.7}
For any $T\in \mathbb T$, we have
\begin{enumerate}
\item If $T$ is the singleton, we have $\psi_T(t)=t$.
\item $\psi_T(0)=0$.
\item
$\psi_T (1)=\begin{cases} 1 \quad & \mbox{if $v(T)=1$}\\
0 \quad & \mbox{if $v(T)\geq 2$}\end{cases}$
\item $\psi_T'(0)=\phi_T$.
\end{enumerate}
\end{lemma}
\pf All statements above follow immediately from (\ref{DefForPsi}),
except the assertion $\psi_T (1)=0$ when $v(T)\geq 2$, which is
exactly (\ref{R00}).  \epfv

\noindent{\bf Forests.} The formula (\ref{DefForPsi}) defines a unique
polynomial $\psi_T(t)$ for each rooted tree $T$.  A {\it forest} is
the disjoint union of finitely many rooted trees.  We extend the
definitions of $\phi_P$ and $\psi_P(t)$ to any forest $P$ as follows:

\begin{defi} \label{forestdef} For any forest $P$ which is the disjoint
union of rooted
trees $T_1,\ldots,T_{k}$, we define $\phi_P$ to be $\phi_{T_1}$ if
$k=1$ and $0$ otherwise.  Define
$\psi_P(t)=\prod_{i=1}^k\psi_{T_i}(t)$.
\end{defi}

\begin{lemma}\label{lemma-3.8}
Let $T$ be a rooted tree with $v(T)\geq 2$. For any proper rooted
subtree
$T'$ of $T$ we have
\BQ
&{}&\psi_{T\backslash T'}(t)=\sum_{k=1}^{v(T)-1} \frac {t^{k}}{k!}
\sum_{\substack {\vec e=(e_1,\ldots,e_{k})\in E(T)^{k}\\
e_1\succ \cdots   \succ e_{k}\\
T_{\vec e,k+1}=T'}}
\phi_{T_{\vec e,1}}\phi_{T_{\vec e,2}}\cdots \phi_{T_{\vec e,k}}
\EQ
\end{lemma}
\pf Let $T^{[j]}$ $(j=1, 2, \dots, d)$ be the connected components of
$T\backslash T'$, and let $e_j^0$ be the edge of $T$ which connects
$T^{[j]}$ with $T'$.  Note that from fixed sequences $e_{j, 1}\succ
e_{j, 2} \succ \cdots, \succ e_{j, k_j}\in E(T^{[j]})$ with
$k_1+k_2+\cdots +k_d=k-d$, appended by the edges $e_j^0$, we can get
$\binom {k}{(k_1+1), (k_2+1),\cdots, (k_d+1)} =\frac
{k!}{(k_1+1)!(k_2+1)!\cdots (k_d+1)!}$ different sequences $e_1\succ
e_2 \succ \cdots \succ e_{k}\in E(T)$ such that $T_{k+1}=T'$.
Therefore, \BQn &{}&\sum_{k=1}^{v(T)-1} \frac {t^{k}}{k!}
\sum_{\substack {\vec e=(e_1,\ldots,e_{k})\in E(T)^{k}\\
e_1\succ \cdots   \succ e_{k}\\
T_{\vec e,k+1}=T'}}
\phi_{T_{\vec e,1}}\phi_{T_{\vec e,2}}\cdots \phi_{T_{\vec e,k}}\nno \\
&=& \sum_{k=1}^{v(T)-1} \frac {t^{k}}{k!}
\sum_{\substack{(k_{1},\ldots,k_{d})\in\bN^{d}\\k_1+k_2+\cdots
+k_d=k-d}}\frac {k!}{(k_1+1)!(k_2+1)!\cdots
(k_d+1)!}\\
&{}& \quad \prod_{j=1}^d
\sum_{\substack{\vec e_{j}=(e_{j,1},\ldots,e_{j, k_j})\in
E(T^{[j]})^{k_{j}}\\
e_{j, 1}\succ \cdots, \succ e_{j, k_j}}}
\phi_{T_{e_{j}, 1}}\phi_{T_{e_{j}, 2}}\cdots \phi_{T_{e_{j}, k_j+1}}\nno
\\
&=&\prod_{j=1}^d \sum_{k_j=0}^{v(T^{[j]})-1}
\frac {t^{k_j+1}}{(k_j+1)!}
\sum_{\substack{\vec e_{j}=(e_{j,1},\ldots,e_{j, k_j})\in
E(T^{[j]})^{k_{j}}\\
e_{j, 1}\succ \cdots \succ e_{j, k_j}}}
\phi_{T_{j, 1}}\phi_{T_{j, 2}}\cdots \phi_{T_{j, k_j+1}}\nno \\
&=&\psi_{T^{[1]}}(t)\psi_{T^{[2]}}(t)\cdots \psi_{T^{[d]}}(t)\,.
\EQn
The last equality follows from (\ref{DefForPsi}).
\epfv

The lemma above allows us to prove the following theorem.  If we let
$\emptyset$ be the empty tree and define $\cP_{\emptyset}=z$, then
Theorem \ref{Theo-3.6} can be seen as the special case $S=\emptyset$
of the theorem below.

\begin{theo}\label{Theo-3.9}
For any rooted tree $S$, we have
\BQ\label{E-3.15}
\exp {(tA)}\cdot{\mathcal P}_S={\cP}_S+
\sum_{T\in \mathbb T} \left ( \sum_{\substack{T'< T \\ T'\simeq S}}
\psi_{T\backslash T'}(t)\right ) \cP_T
\EQ
\end{theo}
\pf
\begin{align}
\exp(tA)\cdot {\mathcal P}_S&=\sum_{k=0}^\infty \frac{t^k}{k!}{A}^k\cdot
{\mathcal P}_S
\nno\\
&=\sum_{k=0}^\infty
\frac{t^k}{k!}\left(\sum_{T\in\bT}\phi_T\cD\right)^k\cdot {\mathcal P}_S\nno
\\
&= {\mathcal P}_S+\sum_{k=1}^\infty
\frac{t^k}{k!}\sum_{(T_1,\ldots,T_k)\in\bT^k}\left[\phi_{T_1}\cD_{T_1}\right]\cdots
\left[\phi_{T_k}\cD_{T_k}\right]\cdot {\mathcal P}_S \nno\\
\intertext{Apply Corollary \ref{key-cor1}, substituting
$Y_{T_i}^{(i)}=\phi_{T_i}$:}
&={\mathcal P}_S+\sum_{k=1}^\infty
\frac{t^k}{k!}\sum_{\substack{T\in\bT\\v(T)\ge1}}\,\,\,\sum_{\substack{\vec{e}
=(e_1,\ldots,e_{k-1})\in
E(T)^{k}\\e_1\succ \cdots \succ e_{k}
\\T_{\vec{e}, k+1}\simeq S}}\phi_{T_{\vec
e,1}}\cdots\phi_{T_{\vec e,k}}\cP_T\nno\\
&={\mathcal P}_S+\sum_{\substack{T\in\bT\\v(T)\ge1}}\,\,
 \left (
\sum_{k=1}^{v(T)-1}
\frac{t^k}{k!}
\,\sum_{\substack{\vec{e}
=(e_1,\ldots,e_{k})\in
E(T)^{k}\\e_1\succ \cdots \succ e_{k}
\\T_{\vec{e}, k+1}\simeq S}}\phi_{T_{\vec
e,1}}\cdots\phi_{T_{\vec e,k}}\right )\cP_T\nno\\
&=\sum_{T\in \mathbb T} \left ( \sum_{\substack{T'\leq T \\ T'\simeq S}}
\psi_{T\backslash T'}(t)\right ){\mathcal P}_T\nno
\end{align}
The last equality follows from Lemma \ref{lemma-3.8}.
\epfv

\begin{propo}\label{prop-3.10}
For any rooted tree $T$, we have \newline

\noindent (a)
\BQ\label{E-3.16}
\psi'_T (t)&=&\phi_T+\sum_{e\in E(T)}
\phi_{T_{e,1}}\psi_{T_{e,2}}(t)
\EQ
(b)
\BQ\label{E-3.17}
\psi'_T (t)
&=&\phi_T+\sum_{S<T} \phi_{S}\,\psi_{T\backslash S}(t)
\EQ
or in other words,
\BQ
\psi'_T (t)&=&\psi_T'(0)+\sum_{e\in E(T)}
\psi_{T_{e,1}}'(0)\,\psi_{T_{e,2}}(t)\nno \\
&=&\psi_T'(0)+\sum_{S<T}\psi'_{S}(0)\,\psi_{T\backslash S}(t)
\EQ
\end{propo}
\pf (a) Applying the chain rule and Theorem \ref{Theo-3.6}, we have
\begin{align}
\frac {\p}{\p t} F_{t}
&= \frac {\p}{\p t} \left(\exp {(tA)}\cdot z\right)\nno\\
&=A\cdot\exp {(tA)}\cdot z \nno\\
&= \left(\sum_{T\in \mathbb T} \phi_T
{\cD}_T\right)\cdot\left(z+\sum_{T\in\mathbb
T}
\psi_T(t){\mathcal P}_T\right)\nno\\
&= \sum_{T\in \mathbb T} \phi_T {\mathcal P}_T+
\sum_{(T_1,T_2)\in \bT^2}  \phi_{T_1}\psi_{T_2}(t)
{\cD}_{T_1}{\cP}_{T_2}\nno \\
\intertext{Applying Corollary \ref{key-cor} with $k=2$, setting
$Y_{T}^{(1)}=\phi_{T}$, $Y_{T}^{(2)}=\psi_{T}(t)$ for all $T\in\bT$:}
&= \sum_{T\in \bT} \phi_T {\mathcal P}_T+
\sum_{T\in \bT} \sum_{e\in E(T)}
\phi_{T_{e,1}}\psi_{T_{e,2}}(t){\cP}_T \nno\\
&=
\sum_{T\in \bT} \left ( \phi_T+\sum_{e\in E(T)}
\phi_{T_{e,1}}\psi_{T_{e,2}}(t)\right ){\mathcal P}_T \nno
\end{align}
But we also have, by Theorem \ref{Theo-3.6},
\BQ\label{E-3.19} \frac {\p}{\p t} F_{t} =\sum_{T\in \mathbb T}
\psi_T'(t){\cP}_T \EQ
Comparing the coefficient of ${\mathcal P}_T$,
and appealing to stable linear independence - specifically, Corollary
\ref{tlic} - we get (\ref{E-3.16}).  (We use the fact that polynomial
functions which agree at all $\alpha\in\bC$ must be equal.)

(b)
\begin{align}
\frac {\p}{\p t} F_{t}
&= \frac {\p}{\p t} \exp {(tA)}\cdot z \nno\\
&= A\cdot\exp(tA)\cdot z \nno\\
&=\exp(tA)\cdot A\cdot z \nno\\
&=\exp(tA)\cdot a \nno\\
\intertext{Applying Theorem \ref{Theo-3.3}:}
&=\exp(tA)\cdot\sum_{S\in\bT}\phi_S\cP_S\nno\\
&=\sum_{S\in\bT} \phi_S\,\exp(tA)\cdot\cP_S\nno\\
\intertext{Applying Theorem \ref{Theo-3.9}:}
&=\sum_{S\in\bT}\phi_S\cP_S+\sum_{S\in\bT}\phi_S
\sum_{T\in\bT} \left(\sum_{\substack{T'<T \\ T'\simeq S}}
\psi_{T\backslash T'}(t)\right)\cP_T\nno\\
&=\sum_{T\in\bT}\phi_T \cP_T+ \sum_{T\in\bT}
 \left(\sum_{S<T}\phi_S
\psi_{T\backslash T'}(t) \right) {\mathcal P}_T\nno\\
&= \sum_{T\in \mathbb T}
 \left ( \phi_T+\sum_{S<T}\phi_S
\psi_{T\backslash T'}(t)\right)\cP_T\nno
\end{align}
Comparing this with (\ref{E-3.19}), and again employing Corollary
\ref{tlic}, we get (\ref{E-3.17}).  \epfv

An interesting consequence of the proposition above is the following
recurrence formula for $\phi_T$ in terms of the number of the leaves of
$T$.

\begin{propo}\label{prop-3.11}
For any rooted tree $T$, suppose that the root $\text{rt}_{T}$
has $d$ children, i.e., $d=\left|\text{rt}_{T}^{+}\right|$.
Then
\BQ\label{E-3.20}
\sum_{r=0}^{l(T)}\,\,\,\,\sum_{\substack{{\{v_1, v_2, \cdots, v_r\}
\subseteq L(T)}\\v_1, v_2, \cdots, v_r\text{ distinct}}}
\phi_{T\backslash \{v_1, v_2, \cdots, v_r\}}=\delta_{d,
1}\phi_{T\backslash \{\text{rt}_{T}\}}
\EQ
\end{propo}
\pf From (\ref{E-3.16}), setting $t=1$, we get
\BQ\label{E-3.21}
\psi'_T(1)=
\begin{cases} \phi_{T}+
\phi_{T\backslash\{\text{rt}_{T}\}} \quad \mbox{if}   \quad d=1\\
 \phi_{T}\quad \quad \quad \mbox{if}   \quad d\geq 2
\end{cases}
\EQ
since, by Lemma \ref{lemma-3.7}, $\psi_{T_2}(1)=0$ except $T_2$ is the
singleton.  From (\ref{E-3.17}), setting $t=1$, we get
\BQ\label{E-3.22}
\psi'_T(1)=\phi_T+\sum_{r=1}^{l(T)}\,\,\,\,\sum_{\substack{{\{v_1, v_2,
\cdots, v_r\}
\subseteq L(T)}\\v_1, v_2, \cdots, v_r\text{ distinct}}}
\phi_{T\backslash \{v_1, v_2, \cdots, v_r\}}
\EQ
since $\psi_{T\backslash S}(1)=0$ except when $T\backslash S$ is the
disjoint union of finitely many
singletons.  Comparing (\ref{E-3.21}) and (\ref{E-3.22}) gives
(\ref{E-3.20}).
\epfv

Before leaving this subsection, we will do some calculations
on the polynomials $\psi_T(t)$ for the
chains $C_n$ and shrubs $S_n$.

Consider the generating functions:
${\cC}(t, x)=\sum_{n=0}^\infty \psi_{C_n}(t)x^n$ (set
$\psi_{C_0}(t)=1$) and
${\cS}(t, x)=\sum_{n=0}^\infty \psi_{S_n}(t)\frac{x^n}{n!}$.

\begin{corol}\label{corol-3.12} The generating functions ${\cC}(t, x)$
and ${\cS}(t, x)$ are given by: \newline

\noindent (a)
\BQ
{\cC}(t,x)=\exp\left(t\ln (1+x)\right)=(1+x)^t
\EQ
 or in other words,
\BQ
\psi_{C_n}(t)=\binom{t}{n}=\frac {t(t-1)\cdots (t-n+1)}{n!}
\EQ
(b)
\BQ
{\cS}(t,x)=\frac {e^{xt}-1}{e^{x}-1}
\EQ
\end{corol}
\pf
(a) By Theorem \ref{Theo-3.6} and Corollary \ref{corol-3.4}. we have
\BQn
{\cC}(t, x)&=&1+\sum_{n=1}^\infty \sum_{k=1}^{v(C_n)} \frac
{t^{k}}{k!}
\sum_{\substack{{\bar e=(e_1,\ldots,e_{k-1})\in E(C_n)^{k-1}}\\
e_1\succ\cdots\succ e_{k-1}}}
\phi_{T_{e,1}}\phi_{T_{e,2}}\cdots \phi_{T_{e,k}}  x^n\\
&=&1+\sum_{n=1}^\infty \sum_{k=1}^{n}\frac{t^{k}}{k!}
\sum_{\substack{(m_1, m_2, \cdots, m_{k})\in\left(\bZ^{+}\right)^{k}\\
m_1+m_2 \cdots+m_{k}=n}}
\frac {(-1)^{m_1}}{m_1}
\frac {(-1)^{m_2}}{m_2}\cdots
\frac {(-1)^{m_{k}}}{m_{k}} x^n\\
&=&   e^{t(-x+\frac {x^2}{2}-\cdots + \frac {(-x)^m}{m}+\cdots)}  \\
&=&\exp\left(t\ln (1+x)\right)
\EQn

(b) Similarly, we have
\begin{align}
{\cS}(t,x)&=\sum_{n=0}^\infty \sum_{k=1}^{v(S_n)}
\frac {t^{k}}{k!}
\sum_{\substack{{\bar e=(e_1,\ldots,e_{k-1})\in E(S_n)^{k-1}}\\
e_1\succ\cdots\succ e_{k-1}}}
\phi_{T_{e,1}}\phi_{T_{e,2}}\cdots \phi_{T_{e,k}}\frac{x^n}{n!}\nno\\
\intertext{Noting that all but one of $\phi_{T_{e,2}}\cdots
\phi_{T_{e,k}}$ are singletons, the remaining one being $S_{n-k+1}$:}
&=\sum_{n=0}^\infty  \sum_{k=1}^{n+1} \frac {t^{k}}{k!}
(k-1)!\binom {m}{k-1} b_{n-k+1}\frac{x^n}{n!}\nno\\
&=x^{-1}\sum_{n=0}^\infty  \sum_{k=1}^{n+1} \frac{(xt)^{k}}{k!}
b_{n-k+1} \frac{x^{n-k+1}}{(n-k+1)!}\nno\\
\intertext{Replacing $n$ by $n-1$:}
&=x^{-1}\sum_{n=1}^\infty  \sum_{k=1}^{n} \frac{(xt)^{k}}{k!}
b_{n-k} \frac{x^{n-k}}{(n-k)!}\nno\\
&=x^{-1}(e^{xt}-1)\frac {x}{e^{x}-1}\nno\\
&=\frac {e^{xt}-1}{e^{x}-1}
\end{align}
\epfv

\begin{rmk}
The formulas of Corollary \ref{corol-3.12} can also be easily derived
from Theorem \ref{Theo-4.2} in the next section.  But we think the
calculations above are more intriguing.
\end{rmk}

\renewcommand{\theequation}{\thesection.\arabic{equation}}
\renewcommand{\therema}{\thesection.\arabic{rema}}
\setcounter{equation}{0}
\setcounter{rema}{0}

\section{The Main Theorem}

\subsection{Main Theorem on $\psi_T(t)$}

In the last section, we defined the polynomial $\psi_T(t)$, for each
rooted tree $T$ (see Theorem \ref{Theo-3.6}).  For each rooted forest
$P$, i.e. the disjoint union of finitely many rooted trees $T_i$
($i=1, 2, \cdots, k$), we also defined $\psi_P$ (see Definition
\ref{forestdef}).  Recalling from \S2.1 the definition of a rooted
subtree, we are now ready to prove the following main
theorem.

\begin{theo}\label{main-theo} Let $t$ and $s$ be indeterminates.
For $T\in\bT$ we have
\BQ\label{additive}
\psi_T(t+s)=\psi_T(t)+\psi_T(s)
+\sum_{T'<T}
\psi_{T\backslash T'}(t)
\psi_{T'} (s)
\EQ
where the last sum runs over all proper rooted subtrees $T'$ of $T$.
\end{theo}
\pf Clearly $\exp((t+s)A)\cdot z
=\exp{(tA)}\cdot\exp{(sA)}\cdot z$, so we have
\begin{align*}
z+\sum_{T\in\bT} \psi_T(t+s){\cP}_T&=\exp((t+s)A)\cdot z=\exp{(tA)}
\cdot\exp{(sA)}\cdot z\\
&=\exp{(tA)}\cdot\left(z+\sum_{T\in\bT}\psi_T(s){\mathcal 
P}_T\right) \\
&=\exp(tA)\cdot z+\exp(tA)\cdot
\left(\sum_{T\in\bT}\psi_T(s){\cP}_T\right)\\
&=z+\sum_{T\in\bT} \psi_T(t){\cP}_T+
\sum_{T\in\bT} \psi_T(s)\left(\exp(tA)\cdot {\cP}_T\right)\\
\intertext{Applying Theorem \ref{Theo-3.9} to $\exp(tA)\cdot {\cP}_T$:}
=z+\sum_{T\in\bT}\psi_T(t)&{\cP}_T
+\sum_{T\in\bT}\psi_T(s)\left({\cP}_T+
\sum_{S\in\bT} \left ( \sum_{\substack{S'< S \\ S'\simeq T}}
\psi_{S\backslash S'}(t)\right)\cP_S\right)\\
=z+\sum_{T\in\bT}\psi_T(t)&{\cP}_T+\sum_{T\in\bT}\psi_T(s){\cP}_T
+\sum_{S\in\bT}\left(\sum_{T\in\bT}\sum_{\substack{S'< S \\ S'\simeq T}}
\psi_{S\backslash S'}(t)\psi_T(s)\right)\cP_S\\
=z+\sum_{T\in\bT}\psi_T(t)&{\cP}_T+\sum_{T\in\bT}\psi_T(s){\cP}_T
+\sum_{S\in\bT}\sum_{S'<S}
\psi_{S\backslash S'}(t)\psi_{S'}(s)\cP_S\\
\intertext{Replacing $S$ by $T$ in the last summation:}
=z+\sum_{T\in\bT}\psi_T(t)&{\cP}_T+\sum_{T\in\bT}\psi_T(s){\cP}_T
+\sum_{T\in\bT}\sum_{T'<T}
\psi_{T\backslash T'}(t)\psi_{T'}(s)\cP_T
\end{align*}
According to Corollary \ref{tlic} and an easy specialization argument,
the theorem follows by comparing the coefficients of ${\mathcal P}_T$ in
the above equation.  \epfv

The {\it difference polynomial} of $g(t)\in\bC[T]$ is defined to be the
polynomial
$\Delta g(t)=g(t+1)-g(t)$.  The following special case of the theorem
above, which
gives a formula for the difference polynomial of $\psi_T (t)$, is most
useful to us.

\begin{theo}\label{Theo-4.2}
For any tree $T$ with $v(T)\geq 2$,  we have
\BQ
\Delta \psi_T (t)&=&\psi_{T_1}(t)\psi_{T_2}(t)\cdots \psi_{T_d}(t)\nno\\
&=&\psi_{T\backslash\{\text{rt}_{T}\}} (t)\nno
\EQ
where $T_i$, $i=1, 2, \cdots, d$ are the connected
components of $T\backslash\{\text{rt}_{T}\}$.
\end{theo}
\pf This follows form Theorem \ref{main-theo} by setting $s=1$ in
(\ref{additive}) and appealing to Lemma \ref{lemma-3.7}, which says
$\psi_T(1)=0$ unless $T$ is the singleton, in which case
$\psi_T(1)=1$.  \epfv

\begin{theo}
For any tree $T$ with $v(T)\geq 2$,  we have \newline

\noindent (a)
\BQ\label{DeltaInLeaves}
\Delta\psi_T (t)=\sum_{r=1}^{l(T)}\,\,\,\sum_{\substack{{\{v_1, v_2,
\cdots, v_r\}
\subseteq L(T)}\\v_1, v_2, \cdots, v_r\text{ distinct}}}\,
\psi_{T\backslash \{v_1, v_2, \cdots, v_r\}}(t)
\EQ
(b)
\BQ
\psi_{T\backslash \{\text{rt}_{T}\}}(t)
=\sum_{r=1}^{l(T)}\,\,\,\sum_{\substack{{\{v_1, v_2, \cdots, v_r\}
\subseteq L(T)}\\v_1, v_2, \cdots, v_r\text{ distinct}}}\,
\psi_{T\backslash \{v_1, v_2, \cdots, v_r\}}(t)
\EQ
where $T_i$, $i=1, 2, \cdots, d$ are the connected
components of $T\backslash\{\text{rt}_{T}\}$.
\end{theo}
\pf
Clearly, (b) follows from (a) and Theorem \ref{Theo-4.2}. For (a),
switch $t$ and $s$ and  set $s=1$
in \ref{additive} to get
\BQn
\psi_T(t+1)=\psi_T(t)+\psi_T(1)+\sum_{T'<T}
\psi_{T\backslash T'}(1)\psi_{T'}(t)
\EQn
By Lemma \ref{lemma-3.7}, we have $\psi_T(1)=0$ and
$\psi_{T\backslash T'}(1)=0$, unless $T\backslash T'$ is a disjoint
union of
singletons, in which case $\psi_{T\backslash T'}(1)=1$. Therefore
\BQn
\psi_T(t+1)-\psi_T(t)=\sum_{r=1}^{l(T)}\,\,\,\sum_{\substack{{\{v_1,
v_2, \cdots, v_r\}
\subseteq L(T)}\\v_1, v_2, \cdots, v_r\text{ distinct}}}\,
\psi_{T\backslash \{v_1, v_2, \cdots, v_r\}}(t)
\EQn
as desired. \epfv

\subsection{Algorithm for $\psi_T(t)$}

>From Theorem \ref{Theo-4.2} we get the following
algorithm for computing $\psi_T(t)$.  Here, for $h(t)\in\bC[t]$,
$\Delta^{-1}h(t)$
is defined to be the unique polynomial $g(t)\in\bC[t]$ such that
$\Delta g(t)=h(t)$ and $g(0)=0$. \vskip 10pt

\noindent {\bf Algorithm.}  For any fixed  rooted tree $T$, we sign a
polynomial $N_v(t)$
to each vertex $v$ of $T$ as follows:
\begin{enumerate}
\item For each leaf $v$ of $T$,  then set $N_v(t)=t$.
\item For any other vertex $v$ of $T$, define $N_v(t)$ inductively
starting
from the highest level by
$N_v(t)=\Delta^{-1} (N_{v_1}(t)N_{v_2}(t)\cdots N_{v_k}(t))$, where
$v_j,$ $j=1, 2, \dots, k$,
 are the distinct
children of $v$.
\end{enumerate}
Then for each vertex $v$ of $T$, $N_v(t)=\psi_{T_v^+}(t)$, where
$T_v^+$ is
the subtree of $T$ rooted at
the vertex $v$. In particular, we have
$\psi_{T}(t)=N_{\text{rt}_{T}}(t)$.
\epfv

The following example applies this algorithm to the shrubs $S_{n}$ to
show that the polynomials $\psi_{S_{n}}(t)$ are closely related to the
Bernoulli polynomials $B_n(t)$ defined by
$\frac{xe^{tx}}{e^{x}-1}=\sum_{n=0}^{\infty}B_{n}(t)\frac{t^{n}}{n!}$.
(Compare this with (b) of Corollary \ref{corol-3.12}.)

\begin{exam}\label{exam-4.4}  Let  $v_1,\ldots,v_n$ be the leaves of
the shrub $S_{n}$.  Following the  algorithm, we first assign the
polynomial $t$ to each leaf $v_{i}$.  The next step in the algorithm
gives
\BQ\label{1.13}
\psi_{S_{n}}(t)=\Delta^{-1} (t^n)
\EQ
One of the fundamental properties of the Bernoulli polynomials
$B_n(t)$ is
\BQ\label{b1}
\Delta B_n(t)=B_{n}(t+1)-B_n(t)=nt^{n-1}\,,
\EQ
and from this and the fact that $\Delta$ commutes with
$\frac{d}{dt}$ one easily derives
\BQ\label{b2}
\frac{d}{dt} B_{n+1}(t)=(n+1)B_{n}(t)\,.
\EQ
>From (\ref{b1}) and (\ref{b2}) we get
\BQ\label{b3}
\Delta^{-1} (t^n)= \int_0^t B_n(u)du =\frac {B_{n+1}(t)}{n+1}-\frac
{B_{n+1}(0)}{n+1}\,.
\EQ
Putting together equations (\ref{1.13}) and (\ref{b3}), we obtain this
relationship between $\psi_{S_{n}}(t)$ and $B_{n+1}(t)$.
\BQ
\psi_{S_{n}}(t)= \int_0^t B_n(u)du =\frac {B_{n+1}(t)}{n+1}-\frac
{B_{n+1}(0)}{n+1}\nno
\EQ
\end{exam}

\subsection{Combinatorial Property of $\psi_T(t)$}

After the main part of this work was done, Professor John Shareshian
pointed out to us that the polynomial $\psi_T(t)$ for rooted trees
coincides with the strict order polynomial $\bar\Omega (P, t)$ for
finite posets (partial ordered sets) $P$ in combinatorics (see
Chapters $3$ and $4$ in \cite{St1}).  We first recall the polynomial
$\bar\Omega (P, t)$ associated with a finite poset, then we show that,
when $P$ is the poset of the set $V(T)$ of vertices of a rooted tree
$T$ with the natural partial order induced by ancestry (the root being
the unique \underbar{smallest} element), we have $\psi_T(t)=\bar\Omega
(P, t)$.

Any rooted tree corresponds in this way to a unique finite poset, and
a finite poset $P$ corresponds to a rooted tree precisely when it
satisfies these two criteria: (1) $P$ has a unique smallest element,
and (2) any interval in $P$ is totally ordered.

For any $n\in \bN$, the chain $C_n$ gives the totally ordered poset
with $n$ elements.  (We can view it as the set $\{1, 2, \cdots, n\}$
with the natural order of the positive integers).  For any poset $P$,
we say a map $f: P\to C_n$ is {\it strict order-preserving} if
$f(a)<f(b)$ in $C_{n}$ whenever $a<b$ in $P$.  It is well-known
that there exists a unique polynomial $\bar\Omega (P, t)$ such that
$\bar\Omega (P, n)$ equals the number of strict order-preserving maps
$f$ from $P$ to $C_n$ for all $n\in\bN$.  This, then, is the theorem
shown to us by John Shareshian.

\begin{theo}\label{Omega}
For any rooted tree $T$, we have
\BQ
\psi_T(t)=\bar\Omega (T, t)
\EQ
(where, on the right, $T$ is viewed as a finite poset as described
above).
\end{theo}
\pf It is obvious that when $T$ is the singleton, $\bar\Omega (T,
t)=t$, hence it is enough to show that the $\bar\Omega (T, t)$ also
satisfies the recursion formula of Theorem \ref{Theo-4.2}.
More precisely, we will show that, in the notation of Theorem
\ref{Theo-4.2}, we have
\BQn
\Delta \bar\Omega (T,n)=\bar\Omega (T, n+1)-\bar\Omega(T,n)
=\bar\Omega(T_1,n) \bar\Omega (T_2,n)
\cdots \bar\Omega(T_d,n)
\EQn
for any $n\in\bN$.

Note that $\Delta \bar\Omega (T,n)$ equals the number of strict
order-preserving maps $f$ from $T$ to $C_{n+1}=\{1,2,\ldots,n+1\}$
such that $f(\text{rt}_{T})=1$.  But this number is also same as the
number of strict order-preserving maps $g$ from $T\backslash
\{\text{rt}_{T}\}$ to $C_{n}$, which is $\bar\Omega(T_1,n) \bar\Omega
(T_2,n) \cdots \bar\Omega(T_d,n)$.  \epfv

\begin{rmk}
It is interesting that the strict order polynomial $\bar\Omega (T,
t)$ for the finite posets induced by rooted trees $T$ can be defined
by a totally different way, namely, according to the formula
(\ref{DefForPsi}) of Theorem \ref{Theo-3.6}.  In fact, this
realization of the strict order polynomial can be generalized to an
arbitrary finite poset $P$.  This generalization and its consequences
will be discussed in the upcoming paper \cite{SWZ}.
\end{rmk}

\renewcommand{\theequation}{\thesection.\arabic{equation}}
\renewcommand{\therema}{\thesection.\arabic{rema}}
\setcounter{equation}{0}
\setcounter{rema}{0}

\section{Some Applications}

For a formal automorphism $F=(F_1, F_2,\ldots,F_n)=z+H$ of
the form identity plus higher, we give a restatement and new
proof of the tree formula for for the formal inverse first proved in
\cite{BCW} and \cite{W}.

\begin{theo}\label{Theo-5.1}
For any rooted tree $T$, we have $\psi_T(-1)=(-1)^{v(T)}$.  Hence the
formal inverse $F^{-1}$ of $F$ is given by\footnote{The formula as
given in \cite{BCW} and \cite{W} did not include the factor
$(-1)^{v(T)}$.  It appears here because of our choice in writing $F=z+H$

instead of $F=z-H$.}
\BQ\label{E-5.1}
F^{-1}=z+\sum_{T\in \bT} (-1)^{v(T)}\,{\cP_{T}}
\EQ
\end{theo}
\pf  The formula (\ref{E-5.1}) follows from $\psi_T(-1)=(-1)^{v(T)}$ by
Proposition \ref{exp} and Theorem \ref{Theo-3.6}.

It is well known in combinatorics (see \cite{St1}) that the strict
order polynomials satisfy $\bar\Omega (T,-1)=(-1)^{v(T)}$, from which
the result follows, in light of Theorem \ref{Omega}.  For
completeness, we give a direct proof here.

We use the mathematical induction on $v(T)$. The case for $v(T)=1$
is trivial.  ($\psi_T(t)=t$ in this case.)  Suppose $v(T)\geq 2$.  By
Theorem \ref{Theo-4.2}, setting $t=-1$, we have
\BQn
\psi_T(0)-\psi_T(-1)=\psi_{T_1}(-1)\psi_{T_2}(-1)\cdots \psi_{T_d}(-1)
\EQn
where $T_i$, $i=1, 2, \cdots, d$ are the connected components of
$T\backslash\{\text{rt}_{T}\}$.
We have $\psi_T(0)=0$, and by induction we may assume the theorem holds
for $T_{1},\ldots,T_{d}$. Hence
\BQn
\psi_T(-1)=-(-1)^{v_{T_1}+v_{T_2}\cdots +v_{T_d}}=(-1)^{v(T)}
\EQn
\epfv

It is known that the Jacobian conjecture (See \cite{BCW} for a
statement of this famous problem.) is equivalent to the assertion
that
\BQ
\sum_{T\in\bT_N} {\cP_{T}}=0
\EQ
for $N>>0$ whenever $H$ is a homogeneous polynomial system (of degree
$\ge2$) and the jacobian determinant $|(D_{j}F_{i})|$ is (everywhere)
non-zero.  In fact, this follows from Theorem \ref{Theo-5.1}, since when

$H$ is homogeneous the polynomials $\sum_{T\in\bT_N} {\cP_{T}}$, for
fixed $N$, are the homogeneous summands of $F^{-1}$ (see Remark
\ref{homdegree}).  When $H$ is homogeneous, the condition
$|(D_{j}F_{i})|=1$ is known to be equivalent to the nilpotence of the
jacobian matrix $JH=(D_{j}H_{i})$ (see \cite{BCW}).  Thus the
following result presents an intriguing statement for comparison.

\begin{propo}\label{prop-5.2}
Assume $H$ is homogeneous of degree $\ge2$.  For any rooted tree, let
$h_{T, k}$ be the number of vertices of height $k$.
Suppose that $(JH)^k=0$.  Then
\BQ\label{E-5.3}
\sum_{T\in \bT_N} h_{T,m} {\mathcal P}_T=0
\EQ
for any $N\in \bN$ and $m\geq k$.
\end{propo}
\pf Suppose that $\deg H=d\geq 2$.  It follows from Euler's formula
that $JH\cdot (z^{\text{t}})=(dH)^{\text{t}}$, from which we get
$\frac1d(JH)^k\cdot(z^{\text{t}})=(JH)^{k-1}\cdot H^{\text{t}}=0$.
(Here the superscript ${}^{\text{t}}$ denotes transpose, converting a
row to a column so that the matrix multiplications make sense.)
For any integer $m\ge1$, a straightforward calculation shows that the
chain $C_{m}$ has the property
${\cP}_{C_{m}}=JH^{m-1}\cdot H^{\text{t}}$.  Therefore
${\cP}_{C_{m}}=0$ for $m\geq k$, and we have
\begin{align*}
0&=\exp {(-A)}\cdot{\cP}_{C_{m}}\\
\intertext{By Theorem \ref{Theo-3.9}, setting $S=C_{m}$ and $t=-1$ in
(\ref{E-3.15}):}
&=
\sum_{T\in \bT}\left(\sum_{\substack{T'\leq T \\ T'\simeq C_{m}}}
\psi_{T\backslash T'}(-1)\right)\,{\cP}_T\\
\intertext{By Theorem \ref{Theo-5.1}:}
&=
\sum_{T\in \bT}\left ( \sum_{\substack{T'\leq T \\ T'\simeq C_{m}}}
(-1)^{v({T\backslash T'})}\right ){\mathcal P}_T\\
&=\sum_{T\in \bT}(-1)^{v(T)-m} h_{T, m} {\mathcal P}_T
\end{align*}
In particular, for any $N\in \bN$, we have
\BQn
\sum_{T\in \bT_N}(-1)^{N-m} h_{T, m} {\mathcal P}_T=
(-1)^{N-m}\sum_{T\in \bT_N} h_{T, m} {\mathcal P}_T=0
\EQn
which gives (\ref{E-5.3}).
\epfv

The proposition above shows that, for a fixed homogeneous polynomial
system $H$, the polynomials ${\cP}_T$ are in some sense quite linearly
dependent to each other.

Finally, let us point out that the formal flow $F_{t}$ gives a formal
flow between $F$ and the identity map $\text{id}\,$, i.e.
$F_{t}|_{{}_{t=1}}=F$ and $F_{t}|_{{}_{t=0}}=\text{id}\,$, having the additional properties $F_t(0)=0$ and $JF_t(0)=I_n$.  It is an open
question in complex analysis whether, for any local analytic map $F$,
such an analytic flow exists.  The usual approach to this question
is to show that $F$ is linearizable, i.e. it is conjugate to a linear
map.  But when $F$ is linearizable the question is still open, even
for the one variable case.  (There are many partial results on this
problem.)  So it is of interest that the formal solution to this
question is given by the very clean formula (\ref{compare}) of Theorem
\ref{Theo-3.6}.  But the question of when $F_{t}$ is locally
convergent is still open.

{\small \sc Department of Mathematics, Washington University in St.
Louis,
St. Louis, MO 63130 }

{\em E-mail}: wright@math.wustl.edu, zhao@math.wustl.edu

\end{document}